\documentstyle[12pt,amssymb]{amsart}
\title[Vector Fields with a nilpotent linear part]
{Integrable analytic Vector Fields With\\ a nilpotent linear part}
\keywords{integrable vector field, nilpotent linear  part, embeddability of 
 mappings}
\subjclass{Primary 58F36, 32S65}
\thanks{Partially supported by NSF grant DMS-9304580 at the Institute for Advanced Study}
\author%[abbreviation of name]
{Xianghong  Gong}
\address{The Institute for Advanced Study\\  School of Mathematics\\ 
Princeton, NJ 08540}
\email{gong@@math.ias.edu}

\newtheorem{thm}{Theorem}[section]
\newtheorem{cor}[thm]{Corollary}
\newtheorem{prop}[thm]{Proposition}
\newtheorem{define}[thm]{Definition}
\newtheorem{lemma}[thm]{Lemma}
\newtheorem{rem}[thm]{Remark}

\newcommand{\be}[1]{\begin{equation}\label{eq:#1}}
\newcommand{\nee}{\end{equation}}
\newcommand{\ben}{\begin{equation*}}
\newcommand{\een}{\end{equation*}}
\newcommand{\nbeq}{\begin{eqnarray}}
\newcommand{\eeq}{\end{eqnarray}}
\newcommand{\bc}[1]{\begin{cor}\label{cor:#1}}
\newcommand{\ec}{\end{cor}}
\newcommand{\bt}[1]{\begin{thm}\label{thm:#1}}
\newcommand{\et}{\end{thm}}
\newcommand{\bl}[1]{\begin{lemma}\label{lemma:#1}}
\newcommand{\el}{\end{lemma}}
\newcommand{\bd}[1]{\begin{define}\label{lemma:#1}}
\newcommand{\ed}{\end{define}}
\newcommand{\bp}[1]{\begin{prop}\label{prop:#1}}
\newcommand{\ep}{\end{prop}}

\newcommand{\an}{analytic}

\newcommand{\co}{coefficient}

\newcommand{\coa}{convergent}

\newcommand{\es}{estimate}
\newcommand{\fo}{formal}
\newcommand{\foa}{formally}
\newcommand{\hol}{holomorphic}

\newcommand{\li}{linearizable}
\newcommand{\seq}{sequence}
\newcommand{\nf}{normal form}
\newcommand{\norm}[2]{\|#1\|_#2}
\newcommand{\ord}{\mbox{ord}\,}
\newcommand{\pows}{power series}

\newcommand{\sinn}{simplification}
\newcommand{\st}{such that}
\newcommand{\tr}{transform}
\newcommand{\trn}{transformation}
\newcommand{\vf}{vector field}

\newcommand{\jq}[1]{\mathopen{<} #1\mathclose{>}}%produce <*,*>
\newcommand{\nd}[1]{\Delta_{(1-#1\theta)r}}
\newcommand{\nnd}[1]{(1-#1\theta)r}
\newcommand{\x}{(x)}
\newcommand{\y}{(y)}
\newcommand{\xy}{(x,y)}
\newcommand{\xyp}{(x^\prime,y^\prime)}
\newcommand{\jz}[1]{\left(\begin{array}#1\end{array}\right)}
\newcommand{\wt}[1]{\mbox{wt}\,#1}
\newcommand{\rl}[1]{Lemma~\ref{lemma:#1}}
\newcommand{\nrc}[1]{Corollary\ref{cor:#1}}
\newcommand{\re}[1]{(\ref{eq:#1})}
\newcommand{\rp}[1]{Proposition~\ref{prop:#1}}
\newcommand{\rt}[1]{Theorem~\ref{thm:#1}}

\begin{document}
%\maketitle
\begin{abstract}
We study the normalization of integrable analytic vector fields with a nilpotent linear part. We prove that such an analytic vector field can be transformed into a certain form by convergent transformations when it has a non-singular 
formal integral. In particular, we show that a formally linearizable analytic vector field with a nilpotent linear part is linearizable by convergent transformations. We then prove that there are smoothly linearizable parabolic analytic 
transformations which cannot be embedded into the flow of any analytic vector field with a nilpotent linear part. 
\end{abstract}
\maketitle
\section{Introduction}
We are concerned with  the normalization of an analytic \vf\ $v$ given by
\be{1}
\frac{dx}{dt}=y+f\xy,\qquad \frac{dy}{dt}=g\xy ,
\end{equation}
where $f,g$ are convergent power series starting with terms of order two.
Since the matrix of  the linear part in \re{1} is nilpotent,
the Poincar\'{e}-Dulac \nf\ gives no \sinn.
In~\cite{takens},
F.~Takens introduced a
\sinn\ for \re{1} by \fo\ \trn s as follows
\be{2}
\frac{dx}{dt}=y+r(x),\qquad \frac{dy}{dt}=s(x).
\end{equation}
The above system is subject to further classification. Using representation theory of 
certain Lie algebras, A.~Baider and J.~C.~Sanders~\cite{baidersanders} gave a complete classification of \re{2}\ 
under suitable non-degeneracy conditions. However, the normal form of
Baider and Sanders excludes the  case  
 $s\x\equiv 0$. We shall see that the vanishing of $s$
 corresponds to the case that $v$ has a non-singular formal integral, i~e.~a 
 \fo\
\pows\ $H\xy $ such that $dH(0)\neq 0$ and $\jq{v,\nabla H}=0$. 

The bifurcation theory on  vector fields in the form \re{1}\ 
was initially studied by F.~Takens~\cite{takens} and
R.~I.~Bogdanov~\cite{bogdanov}. We also refer
to~\cite{arnoldilyashenko} for the survey on \vf s with a
nilpotent linear part.
In this paper we shall deal with the convergence of the initial normalization
given by Takens for integrable \vf s. We shall  prove the following.
\bt{1} 
Let $v$ be an
analytic \vf\ given by $\re{1}.$ Assume that $v$ has a non-singular formal
integral$.$ Then $v$ can be \tr ed into $\re{2}$ with
$s\equiv 0$ through a  convergent \trn$.$
%\be{3}
%\frac{dx}{dt}=-y+h\x,\qquad \frac{dy}{dt}=0,
%\end{equation}
%Where $h\x\equiv 0$, or 
%\be{4}
%h\x= \epsilon x^\sigma+\dots,\qquad \epsilon=1 \ \mbox{for}\ 
%\sigma\ \mbox{even, and}\
% \epsilon=\pm 1\ \mbox{for}\ \sigma\ \mbox{odd}.
%\end{equation}
\et
We cannot give a complete classification for \re{1} under the assumption
of formal integrability. As a partial normalization,
 we have
\bt{2}
Let $v$ be an \an\ \vf\ given as in $\rt{1}.$ Assume  that $f,g$ in
$\re{1}$ are real power series$.$
 Then
there is a real formal \trn\ which
transforms $\re{1}$ into
\be{4}
\frac{dx}{dt}=y+r^*(x),\qquad
%\epsilon x^\sigma+\sum_{j>\sigma, \sigma\nmid j}b_jx^j,\qquad
\frac{dy}{dt}=0,
\end{equation}
where $r^*\equiv 0,$ or
$$
r^*(x)=\epsilon x^\sigma+\sum_{j>\sigma, \sigma\nmid  j}r_j^*x^j,
$$
and $\epsilon=1$ for $\sigma$ even$,$ $\epsilon=\pm 1$ for $\sigma$ odd$.$
Furthermore,  when $\sigma$ is even, all
$r_j^*\ (\sigma<j<2\sigma)$ are invariants$;$
when $\sigma$ is odd$,$ $\epsilon$ and  the coefficients
$r_j^*\ (\sigma<j<2\sigma)$ are invariants$,$ of which
the first non-zero coefficient of even order is normalized to be positive$.$
\et
We shall see that there are infinitely many invariants for \re{4}\ if
$2<\sigma<\infty$.
We are unable to determine whether  \re{4}\ can be realized by convergent transformations.
Neither shall we  deal with the convergence of the normal form of Baider-Sanders. In fact, the convergence
proof for the vector fields considered here  depends essentially on the
assumption of the existence of non-singular integrals.

\rt{1}
demonstrates a significant
 difference  between   real analytic \vf s with a nilpotent linear part and
real analytic mappings with a  unipotent linear
part, i.~e. parabolic mappings. The parabolic mappings arise naturally
from  glancing hypersurfaces considered by R.~B.~Melrose~\cite{melrose1},
and also
 real analytic Lagrangian surfaces with a complex tangent  studied by S.~M.~Webster~\cite{websterpisa}.
 In~\cite{gong1}, it was shown that there exist 
 real analytic \trn s which are \foa\ equivalent to the linear parabolic
mapping $T\xy =(x+y, y)$, but they are not \li\ through any \coa\ \trn. On
the other hand, \rt{1} shows
that  a \vf\ \re{1}\ is \li\ by convergent \trn s if
and only if it is formally \li. To  further understand the
 distinct nature of  normalizing \an\ \vf s with nilpotent linear parts and parabolic
\an\  mappings, we shall prove  the following result.
\bc{3}
There exists a  smoothly \li\ real analytic \trn\ $\varphi=T+O(2)$ which is not
 embeddable in any
 neighborhood of the origin as 
the time-1 mapping of any  real analytic \vf\ $\re{1}.$
\ec

This paper is organized as follows. In section 2, we shall consider the formal theory of the integrable vector field \re{1}. We shall also characterize the integrability in terms of the formal normal form as well as the singular points of the vector fields. The proof of \rt{2} will be given in section 2.
 In section 3, we shall prove the convergence of  solutions to the approximate equations  arising from the normalizing the vector fields in \rt{1}. We complete the proof of \rt{1} in section 4 through a KAM argument. In section 5, we shall discuss the embeddability of a parabolic mapping into the flow of a vector field \re{1}. The proof of \nrc{3} will be presented 
in section 5.

We would like to thank the referee for bringing the article of D.~Cerveau and R.~Moussu~\cite{cer} to our attention. In terms of the vector fields, the results in~\cite{cer} contain 
  a complete holomorphic (convergent) classification for  all  ideals of holomorphic vector fields  in ${\Bbb C}^2$ which  are formally equivalent to a given ideal of vector fields 
defined by $\omega=0$,  where the holomorphic $1$-form $\omega$ is a  certain perturbation of  $ydy-x^ndx$.

\section{\fo\ normalizations}\setcounter{equation}{0}
\label{section:2}
In this section, we shall first construct a  \fo\ \trn\ $\Phi$  which \tr s \re{1}\ into \re{2}.  
The convergence of $\Phi$ will be determined in section 4.  We shall
 also discuss
the formal integrability. Finally, we shall give a proof for \rt{2}.

Throughout the discussion of this paper, we shall decompose  a \pows\ $p\xy $ into the following form
$$
p\xy =p_0\x+p_1\xy ,\qquad p_0\x=p(x,0).
$$
Then the mapping $p\to p_1$ defines a projection $\Pi\colon {\Bbb C}[[x,y]]\to y{\Bbb C}[[x,y]]$, where
${\Bbb C}[[x,y]]$ is the ring of formal power series in $x,y$, and $y{\Bbb C}[[x,y]]$ is the ideal generated by $y$. For a non-zero power series $p\xy$, we denote by $\ord p$ the largest
integer $k$ such that $p_{\alpha,\beta}=0$ for $\alpha+\beta<k$.
We put $\ord {p}=\infty$ when $p\equiv 0$. We also write
$$
p(x,y)=O(k),
$$
if $ \ord{p}\geq k.$

Let $\varphi$ be a \trn\ defined by
\be{2.1}
x^\prime=x+u\xy ,\qquad y^\prime=y+v\xy ,
\end{equation}
where $u,v$ are \pows\ starting with terms of order two. We say that $\varphi$ is {\it normalized\/} if
\be{2.2}
u(0,y)\equiv 0\equiv v(0,y).
\end{equation}
One can see that $\varphi$  is a normalized \trn\ if and only if  $\varphi$ preserves the 
$y$-axis, and its restriction to the $y$-axis is the identity mapping. Therefore, the normalized \trn s form a group.

We shall seek a unique normalized \trn\ 
$$
\Phi\colon  x^\prime=x+U\xy ,\qquad y^\prime=y+V\xy ,
$$ 
which transforms \re{1}\ into \re{2}.  This leads to  the following functional equations:
\begin{eqnarray}
yU_x\xy  -r\x&=&V\xy -f \xy +E_1\xy ,\label{eq:2.3}\\
yV_x\xy  -s\x&=& -g \xy +E_2\xy ,\label{eq:2.4}
\end{eqnarray}
where
%\be{2.5}\begin{split}%{array}{l}%
\begin{eqnarray}
E_1\xy& =r(x+U )-r\x-f U_x\xy -g U_y\xy,\nonumber%\label{eq:2.5}
\\
E_2\xy& =s(x+U )-s\x-f V_x\xy -g V_y\xy.\label{eq:2.5+}
%\end{split}%{array}
%\end{equation}%
\end{eqnarray}

We shall prove that under the normalizing condition
\re{2.2},  the equations \re{2.3} and \re{2.4} have a   unique   solution $\{U,V,r,s\}$.  Let us denote
by $E_{j;\alpha,\beta}$ the \co\ of the term $x^\alpha y^\beta$ of $E_j$.  Then, it is easy to see that for given $\alpha+\beta=n$,  $E_{j;\alpha,\beta}$  is a polynomial in $r_{\alpha^\prime}, s_{\alpha^\prime}\ (\alpha^\prime<n)$ and $f_{\alpha^\prime,\beta^\prime}, g_{\alpha^\prime,\beta^\prime}\ (\alpha^\prime+\beta^\prime<n)$ with integer \co s.
Comparing the \co\ of $x^{\alpha-1} y^{\beta+1}$ on the both sides of  \re{2.4}, we get
\be{v}
V_{\alpha,\beta}=\frac{1}{\alpha}\left(-g_{\alpha-1,\beta+1}+E_{2;\alpha-1,\beta+1}\right),
\quad  1\leq \alpha\leq n.
\end{equation}
By the normalizing condition \re{2.2}, we also have 
\be{0}
U_{0,n}=0,\quad V_{0,n}=0.
\end{equation}
Next, we compare the \co\ of $x^{\alpha-1} y^{\beta+1}$ on the both sides of  \re{2.3},  and obtain
\be{u}
U_{\alpha,\beta}=\frac{1}{\alpha}\left(
V_{\alpha-1,\beta+1}-f_{\alpha-1,\beta+1}+E_{1;\alpha-1,\beta+1}\right),\quad  1\leq \alpha\leq n.
\end{equation}
Finally,  the \co of $x^n$ on the both sides of \re{2.3} and \re{2.4} gives us
\be{rs}
r_n=f_{n,0}-E_{1;n,0}-V_{n,0},\quad s_n=g_{n,0}-E_{2;n,0}.
\end{equation}
Notice that $E_{j;\alpha,\beta}=0$ for $\alpha+\beta=2$. Thus, the \co s of  $r,s,U,V$  of order $2$ are uniquely determined by the \co s of $f,g$ of order $2$ through formulas 
\re{v}-\re{rs}.  By induction, one can show that the \co s of $r,s,U,V$ of order $n$ are uniquely 
determined by   the \co s of $f,g$ of order $n$. Therefore, there exists a unique solution 
$\{r, s,U,V\}$ to \re{2.3} and \re{2.4}, of which  $U,V$ satisfy the condition \re{2.2}. This proves that
\re{1} can be transformed into \re{2} by
 a unique normalized formal transformation.

 For the late use, we  observe that  if \re{2.3} and
\re{2.4}\  are solvable for  $r,s,U,V$ with $s\equiv 0$, then  $g$ must satisfy the condition
\be{late}
\ord  g_0\geq\ord  g_1.
\end{equation}
To see this, we put $s\equiv 0$ in   \re{2.4} and (\ref{eq:2.5+}).  From (\ref{eq:2.5+}), we 
see that $\ord E_2\geq \ord V+1$. 
By \re{2.2}, we have $\ord V=\ord \{yV_x\xy\}$. Applying $\Pi$ to \re{2.4}, we then get $\ord V\geq
\ord g_1$.  By comparing the orders on both sides of \re{2.4}, one can see 
easily that
 \re{late} holds.

Next, we want to describe the integrability. We have
\bp{2.1}
Let $v$ be an analytic \vf\  defined by $\re{1}.$ Assume that $\re{1}$ is transformed into  $\re{2}$ through a \fo\ \trn\ $\varphi.$ Then the following are equivalent$:$
\newcounter{eq}
\begin{list}
{$(\alph{eq})$}{\usecounter{eq}}
\item 
 $s(x)\equiv 0.$
\item $v$ has a non-singular formal integral$.$
\item  The set of singular points of $v$ is a curve through the origin$.$
\end{list}
\ep
\begin{rem} There exist integrable \vf s with a nilpotent linear part,
which have an isolated singular point at the origin. For instance, consider
$$
\frac{dx}{dt}=y,\qquad\frac{dy}{dt}=x^2.
$$
Then $H\xy=2x^3-3y^2$ is an integral of the system.
Hence, it is essential that the integral in \rt{1}\ is 
non-singular.
 \end{rem}
From \rt{1} and \rp{2.1}, we have the following.
\bc{+}
Let $v$ be an analytic vector field given by $\re{1}.$ Then $v$ is linearizable by convergent transformations if and only if it is formally linearizable.
\ec

\begin{pf*}{Proof of \it{\rp{2.1}}} It is obvious that (a)$\Rightarrow$(b). 
  To show that (b)$\Rightarrow$(a), we assume that $s(x)$ does  not vanish identically and write
$$
s(x)=s_\tau x^\tau+\ldots,\qquad s_\tau\neq 0.
$$
 We want to prove that $\tau$ is the largest positive integer $k$ such that
there exists a  formal power series $H$ with
\be{ch}
\ord {(\jq{v, \nabla H})}\geq k,\qquad dH(0)\neq 0.
\end{equation}

First, the existence of such power series $H$ does not depend on the
choice of formal
coordinates for $v$. Hence, we may assume that $v$ is a formal \vf\ in the form \re{2}.
Moreover, if  $H\xy=y$, then
$\ord {(\jq{v, \nabla H})}=\tau$. 

Next, we want to show that there is no power  series $H$ such that \re{ch}\
holds for some $k>\tau$. Assume that such a power series $H$ exists, and put 
$$
H\xy=\sum_{j=0}^\infty H_j\y x^j.
$$
From \re{ch}, we get
$$
yH_x\xy=-r\x H_x\xy-s\x H_y\xy+O(k).
$$
Expanding both sides as power series in $x$ and comparing the \co of $x^{j-1}$, we obtain
\be{h_n}
jyH_j\y=-\sum_{l=2}^{j-1}(j-l)r_lH_{j-l}\y
-\sum_{l=\tau}^{j-1} s_lH_{j-l}^\prime\y+O(k-j+1)
\end{equation}
for  $1\leq j\leq k$. In particular, we have 
\be{h_n+}
H_j\y=O(k-j) 
\end{equation}
for $j=1,2$. In fact, we want to show that \re{h_n+} holds for $1\leq j\leq \tau$. For the induction, 
we assume that \re{h_n+} holds for $1\leq j\leq\tau^\prime, \tau^\prime<\tau$. Then for $2\leq l\leq
\tau^\prime$, we get $H_{\tau^\prime-l}=O(k-\tau^\prime-1+l)$. It is clear that $k-\tau^\prime-1+l
\geq k-\tau^\prime+1$ for $l\geq 2$. Hence
$$
H_{\tau^\prime-l}=O(k-\tau^\prime+1).
$$
 Set $j=\tau^\prime+1$ in \re{h_n}. Then   the first summation on the right side 
of \re{h_n} can be replaced by $O(k-\tau^\prime+1)$. The second summation in \re{h_n} vanishes, 
since $\tau^\prime<\tau$.  Hence, \re{h_n+}  holds for $j=\tau^\prime$. Therefore, we have verified by induction  that \re{h_n+} holds for $1\leq j\leq \tau$. We now take $j=\tau+1$ for \re{h_n}, and  get
$$
(\tau+1)yH_{\tau+1}\y=-s_\tau H_0^\prime\y+O(k-\tau).
$$
Since $dH(0)\neq 0$, we have $H_0^\prime(0)\neq 0$. Hence, the above identity
cannot hold if $k>\tau$. Therefore,  the order $\tau$ of $s$ is an invariant. 

Now, the formal integrability implies that there  is a formal power series
$H$ such that \re{ch}\ holds for all integer $k$, which yields
$s\equiv 0$.

(b)$\Leftrightarrow$(c). Notice that the origin is not an isolated singular point of $v$,  if and only if $g\xy $ can be divided by $y+f\xy $. We have
\be{jacob}
\left(\begin{array}{c} y^\prime+r(x^\prime)\\ s(x^\prime)\end{array}\right)=D\varphi\xy \cdot
\left(\begin{array}{c} y+f\xy \\ g\xy \end{array}\right),
\end{equation}
where $D\varphi$ is the Jacobian matrix  of the transformation $(x^\prime,y^\prime)=\varphi(x,y)$.  We first assume that (b) holds. Then \re{jacob} implies that $y^\prime+r(x^\prime)$ and $s(x^\prime)$ can be divided by $y+f\xy$ for $\xy=\varphi^{-1}(x^\prime, y^\prime)$. Since $y^\prime+r(x^\prime)$ is irreducible, then $s(x^\prime)$ must be divided by $y^\prime+r(x^\prime)$. 
Hence, we can write
\be{sa}
s(x^\prime)=a(x^\prime, y^\prime)(y^\prime+r(x^\prime))
\end{equation}
for some convergent power series $a(x^\prime, y^\prime)$. Assume for contradiction that $a\neq 0$. 
As power series in $x^\prime, y^\prime$, we compare orders on both sides of \re{sa} and get
 $$
\ord s=\ord a+1.
$$
Next, we set $y=0$ in \re{sa}. Then as power series in $x$ along, we obtain
$$
\ord s\geq\ord a+2.
$$
Thus, the contradiction implies that $a\equiv 0$, i.~e. $s\equiv 0$. Conversely, let us assume that 
(c) holds.  Then from \re{jacob}, we see that both $y+f\xy$ and $s(x)$ can be divided by $s\circ\varphi\xy$. Since $y+f\xy$ is irreducible, then $s(x)$ must
 be divided by $y+f\xy$. Hence, we get (b).  This completes the proof of \rp{2.1}. 
\end{pf*}

Next, we assume that  $s(x)\equiv 0$.
By a linear \trn
$$
\xy\to(ax,ay),
$$
one can achieve that
$$
r_\sigma=\epsilon=\left\{\begin{array}{ll}
1,& \mbox{if $\sigma $ is even},\vspace{.5ex}\\
\pm 1,& \mbox{if $\sigma $ is odd}.
\end{array}\right.
$$
With the normalization for $r_\sigma$, we have $a=1$ if $\sigma$ is even, and $a=\pm 1$ if $\sigma$ is odd.
 
From the proof of \rp{2.1}, we see that the order of vanishing of $s(x)$ is an
invariant for the system \re{2}. 
One can also give a characterization for the order of vanishing of 
 $r(x)$ when $s\equiv 0$. Here, we need \rt{1}. First, \rt{1}
implies that $v$ is actually integrable.
The curve of the singular
points of $v$ is
$$
S\colon y+f\xy=0.
$$
Let $\gamma$ be the level curve of integral passing through the origin.
Then $\sigma $ is the order of  contact of $\gamma$ and $S$ at the origin. When
$\sigma$ is even, $S$ is located on one  side of $\gamma$. If
$\sigma$ is odd and $\epsilon=-1$, then the orbits of the \vf\ are
attracted to $S$. When $\sigma$ is odd and $\epsilon=+1$,  the
orbits will leave $S$ along level curves of integral.

%We now turn to the proof of \rt{2}. 

%\vspace{1.5ex}\noindent
\begin{pf*}{Proof of \it{\rt{2}}} We put $i+j\sigma$ to be the {\it weight } of $x^iy^j$.
For a power series $p\xy$, we denote  by $p_{n;j}$ the \co\ of $x^iy^j$
with weight $n$, and also by $\wt{p}$ the largest integer $n$ such that all
coefficients of $p$ of weight less than $n$ vanish.

Assume that $\varphi$ is  a \trn\ which \tr s \re{2}\ into 
$$
\frac{d\xi}{dt}=\eta+\epsilon \xi^\sigma+\sum_{j=\sigma+1}^\infty r_j^*\xi^j,\qquad \frac{d\eta}{dt}=0.
$$
Then $\varphi$ has the form
\be{6.0}
\left\{\begin{array}{l}
\xi=ax+u\xy,\qquad \wt{u}\geq 2,\vspace{.5ex}\\
\eta=ay+v(y), \end{array}\right.
\end{equation}
where $a=1$ if $\sigma$ is even, and $a=\pm 1$ if $\sigma$ is odd. Notice that $a^{\sigma-1}=1$.
Then we have the following functional equation
\be{6.1}
yu_x+\epsilon x^\sigma u_x-\epsilon \sigma x^{\sigma-1}u-r_0^*(ax)+ar_0(x)=v\y+E\xy,
\end{equation}
where
\begin{gather*}r_0\x=\sum_{j>\sigma}r_jx^j,\qquad
r_0^*\x=\sum_{j>\sigma}r_j^*x^j,\\
E\xy=r_0^*(ax+u)-r_0^*(ax)-r_0u_x+\epsilon \left((ax+u)^\sigma-ax^\sigma-
\sigma x^{\sigma-1}u\right).
\end{gather*}

We want to show that $\wt{u}\geq \sigma$. Let
$$u\xy=u_{k;0}x^k+\ldots,\qquad 1<k<\sigma.
$$
Notice that the weight of terms in $E$ is at least $k+\sigma$. On the other
hand, the coefficient of $x^{k-1}y$ on the left side of \re{6.1} is
$ku_{k;0}$.  Hence, $u_{k;0}=0$.  This shows that the weight of $u$ is at
least $\sigma$. By collecting terms in \re{6.1} with weight less than $2\sigma-1$, we get
$$
r_j^*a^j=a r_j,\qquad \sigma<j<2\sigma-1.
$$
Next, by comparing the coefficients of  terms of weight $2\sigma-1$ in \re{6.1}, we obtain
 $$
u_{\sigma;0}=\epsilon u_{\sigma;1},\qquad r_{2\sigma-1}^*=r_{2\sigma-1}.
$$
Therefore, the coefficients $r_j\ (\sigma+1\leq j\leq 2\sigma-1)$ are invariants, if we restrict the first non-zero coefficient of even order to be positive when $\sigma$ is odd.   

We now assume that $a=1$ in \re{6.0}. To achieve that  $r_{(m+1)\sigma}^*=0$ for $m\geq 1$, we
compare the coefficients of weight $(m+1)\sigma$ on both sides of \re{6.1},
which gives
\begin{gather*}
u_{n;0}=\dfrac{\epsilon}{\sigma-n}
\left(r_{n+\sigma-1}-E_{n+\sigma-1;0}\right),\quad n=m\sigma+1,\\
 u_{n;j}=\dfrac{E_{n+\sigma-1;j}-\left(n-(j-1)\sigma\right)u_{n;j-1}}
{\epsilon\left(n-(j+1)\sigma\right)},\quad 1\leq j\leq m,\\
v_{m+1}=u_{n;m}-E_{n+\sigma-1;m+1}.
\end{gather*}
Therefore,  the coefficients of $u$ of weight $m\sigma+1$ and the coefficient of $v$ of weight $(m+1)\sigma$ are uniquely determined by $r_k\ (k\leq(m+1)\sigma)$ and the coefficients of $u$ with weight less than $m\sigma+1$. 

Next,  consider the coefficients of $u$ of  weight 
$$n=m\sigma+k,\qquad   0\leq k< \sigma,\quad k\neq 1.
$$
 From \re{6.1}, we get 
\begin{gather*}
ku_{n;m}=E_{n+\sigma-1;m+1},\qquad k\geq 1,\\
u_{n;j-1}=\dfrac{E_{n+\sigma-1;j}-\epsilon\left(n-(j+1)\sigma\right)u_{n;j}}{n-(j-1)\sigma},\quad 1\leq j\leq m, \\
r_{n+\sigma-1}^*=r_{n+\sigma-1}+\epsilon\left(n-\sigma\right)u_{n;0}-E_{n+\sigma-1;0}.
\end{gather*}
Hence, $r_{(m+1)\sigma+k-1}^*$ and  coefficients of $u$  with weight
$m\sigma+k$, except for $u_{m\sigma;m}$,
are uniquely determined by $r_j$ with $j\leq
(m+1)\sigma+k-1$ and the
coefficients of $u$ with weight less than $m\sigma+k$.
Therefore, one can achieve that $r_{j\sigma}^*=0$ for $j=2,3,\ldots,m$ through a formal
transformation \re{6.0}, of which the \co s of $u$ with weight
up to  $ (m-1)\sigma+1$ and \co s of $v$ with weight  up to $m\sigma$ are
uniquely determined by  the coefficients $u_{\sigma;1},u_{2\sigma;2},\ldots,
u_{(m-1)\sigma;m-1}$. 
% The proof of \rt{2}\ is complete. 
\end{pf*}

Furthermore, by counting the number of coefficients, we see that the system \re{4}\ has  infinitely many invariants when
$\sigma\geq 3$.

\section{Solutions to  approximate equations}\setcounter{equation}{0}
\label{section:3}
The convergence of \trn\ $\Phi$ cannot be determined directly from the functional 
equations \re{2.3}\ and \re{2.4}.  In this section, we shall  give some \es s  of  solutions to the approximate equations.

We shall consider the following 
 approximate equations 
\begin{eqnarray}
-yu_x\xy +\widetilde{f}_0\x &=&f\xy +f_0\x u_x\xy
\label{eq:3.1}\\
& & -\ f_0^\prime\x u\xy-v\xy ,\nonumber\\
-yv_x\xy +\widetilde{g}_0\x &=&g\xy +f_0\x v_x\xy , 
    \label{eq:3.2}
\end{eqnarray}
in which $\tilde{f}_0$ and $\tilde{ g }_0$ are added to adjust terms purely in $x$.

One can see that under the normalizing condition \re{2.2}, power series
$u,v,\tilde{f}_0,\tilde{g}_0$ are determined uniquely from \re{3.1}\ and \re{3.2}.   The proof can be
given by an  argument similar to the proof of the existence and the uniqueness of solutions $r,s,U,V$ to \re{2.3} 
and \re{2.4}. We left the details to the reader. 

With the above solution $\{u,v\}$, we define a formal transformation   $\varphi$ by \re{2.1}. Assume that $\varphi$  \tr s \re{1}\  into 
\be{n1}
\frac{dx^\prime}{dt}=y^\prime+p(x^\prime, y^\prime),\qquad 
\frac{dy^\prime}{dt}=q(x^\prime,y^\prime).
\end{equation}
Then we have the following identities:
\begin{eqnarray}
p(x^\prime,y^\prime)& =&f\xy  +(y+f\xy)u_x\xy  +g\xy u_y\xy -v\xy,\label{eq:3.12}\\
 q(x^\prime,y^\prime)&  
=& g\xy +(y+f\xy )v_x\xy +g\xy v_y\xy \label{eq:3.13}
\end{eqnarray}
with $\xy=\varphi^{-1}(x^\prime,y^\prime) $. We denote
$$
d_0=\min\{\ord  f_1,\ord 
g\}, \quad d_1=\min\{\ord  p_1,\ord  q\}.
$$

 We  need the following.
\bl{4.1}
Assume that $f,g$ in $\re{1}$ are holomorphic in $\Delta_r.$ 
Let $u,v$ be solutions to $\re{3.1}$ and $\re{3.2},$ which satisfy the
normalizing  condition $\re{2.2}.$ Assume that $\varphi$ defined by $\re{2.1}$ 
\tr s $\re{1}$ into $\re{n1}.$ 
 If $v$ defined by $\re{1}$ has a non-singular formal integral$,$  then
\be{4.1}
d_1\geq 2d_0-1.
\end{equation}
\el
\begin{pf}
 Applying $\Pi$ to \re{3.2}, we get
$$
\ord  \{yv_x\}\geq \min\{\ord g, \ord \{f_0v_x\}\}.
$$
 By \re{2.2}, we also have $\ord v_x=\ord v-1$. Notice that $\ord f\geq 2$. Hence 
\be{ov}
\ord v\geq d_0.
\end{equation}
 We now apply $\Pi$ to \re{3.1} and get 
\be{ou}
\ord u\geq d_0. 
\end{equation}
 Adding  $yu_x\xy-f_0(x)$ to both sides of \re{3.1}, we see that
\be{3.3}
  \ord  (\tilde{f}_0-f_0) \geq d_0.
\end{equation}

Eliminating $v$ from \re{3.1}\ and (\ref{eq:3.12}), we get
$$
p(x+u ,y+v )=\tilde{f}_0(x)+f_0^\prime u\xy +f_1  u_x\xy +g  u_y\xy .
$$
From \re{ov} and \re{ou}, we see that
\begin{gather}
p(x,y)=\tilde{f}_0(x)+O(d_0),\label{eq:4.2}\\
p(x+u ,y+v )=\tilde{f}_0(x)+f_0^\prime(x)u\xy +O(2d_0-1).\label{eq:4.3}
\end{gather}
Combining \re{4.2}\  and \re{3.3}, we obtain that $p\xy =f_0\x+O(d_0)$, i.~e.
$$p_1\xy =O(d_0),\qquad p_0\x=f_0\x+O(d_0).$$
Thus, we get
$$p(x+u ,y+v )=p\xy +f_0^\prime\x u\xy +O(2d_0-1).$$
Now \re{4.3}\ yields
 $$
p\xy =\tilde{f}_0\x+O(2d_0-1).
 $$
In particular, we have
\be{4.4}
\ord {p_1}\geq 2d_0-1.
\end{equation}
From \re{3.2}\ and \re{3.13}, it follows that
$$q(x+u ,y+v )=\tilde{g}_0(x)+O(2d_0-1).$$
This implies that $\ord  q\geq d_0$, and 
$$q(x+u ,y+v )=q\xy +O(2d_0-1).$$
 Therefore, we have
$q\xy =\tilde{g}_0(x)+O(2d_0-1).$
In particular, we obtain that
\be{4.5}
\ord  q_1\geq 2d_0-1.
\end{equation}

Notice that \re{n1}\ also has a non-singular formal integral. From \re{late}, it
follows that
$$
\ord  q_0\geq
\ord  q_1.
$$
Hence, \re{4.4}\ and \re{4.5}\ give us \re{4.1}.
% The proof of \rl{4.1}\ is complete.
\end{pf}

We now want to show the convergence of $u$ and $v$. 
We need the following  notations
%\begin{gather*}
%a_n=|f_{0,n+1}|,\qquad 
%b_n=\max\{|f_{\alpha,\beta}|,|g_{\alpha,\beta}|;
% \alpha+\beta=n,\alpha\geq 1\},\\
%\mu_n=\max\{|u_{\alpha,\beta}|;
% \alpha+\beta=n,\alpha\geq 1\},\quad
%\nu_n=\max\{|v_{\alpha,\beta}|;
% \alpha+\beta=n,\alpha\geq 1\}.
%\end{gather*}
%\begin{alignat*}2
\begin{equation*}
%$$
\begin{array}{ll}
a_n =|f_{n+1,0}|,&\quad \mu_n=\displaystyle{\max_{\alpha+\beta=n,\alpha\geq
1}\{|u_{\alpha,\beta}| \}},\vspace{.75ex}
\\ %\vspace{1.5ex}
b_n  =\displaystyle{\max_{\alpha+\beta=n,\beta\geq 1}}\{|f_{\alpha,\beta}|,|g_{\alpha,\beta}|\},&\quad
\nu_n=\displaystyle{\max_{\alpha+\beta=n,\alpha \geq 1}\{|v_{\alpha,\beta}|
\}}.
\end{array}%$$
\end{equation*}
%\end{alignat*}
  Given two power series $p\xy $ and $q\xy $,  we shall denote 
$
p\prec q,$ if $ |p_{\alpha,\beta}|\leq q_{\alpha,\beta}$ for all $\alpha, \beta\geq 0$. 
We shall also denote
$$
\hat p\xy =\sum |p_{\alpha,\beta}|x^\alpha y^\beta.
$$
 
 Comparing the \co\ of $x^\alpha y^\beta$ on both sides of  \re{3.2}, we have
$$
-(\alpha+1)v_{\alpha+1,\beta-1}=g_{\alpha,\beta}+
\sum_{\alpha^\prime+\alpha^{\prime\prime}= \alpha}(\alpha^\prime+1)v_{\alpha^\prime +1,\beta}f_{\alpha^{\prime\prime},0}    
$$
 for  $\beta\geq 1$ and $\alpha+\beta=n$. Let $\gamma=\alpha^{\prime\prime}
-1$. Then $\alpha^\prime+1+\beta=n-\gamma$. Hence
$$
|v_{\alpha+1,\beta-1}|\leq b_n+\sum \nu_{n-\gamma}a_\gamma.
$$
Therefore, we have
$$
\nu(t)\prec a(t)\nu(t)+b(t),
$$
which gives us
\be{3.4}
\nu(t)\prec \frac{b(t)}{1-a(t)}.
\end{equation}

Solving  \re{3.1}\ for $u$,  one gets
$$
\mu(t)\prec \nu(t)+2a(t)\mu(t)+b(t).
$$
From  \re{3.4}, it follows that
$$
\mu (t)\prec \frac{1}{1-2a(t)}\left(\frac{b(t)}{1-a(t)}+b(t)\right),
$$
which yields
\be{3.5}
\mu(t)\prec\frac{2b(t)}{\left(1-2a(t)\right)^2}.
\end{equation}

\begin{rem} The results in this paper, except for \rt{2}\ which needs  obvious modifications, 
 are valid for holomorphic \vf s. In fact, for the proof of \rt{1}, we shall
introduce holomorphic coordinates. \end{rem}

 From now on, we shall treat all
variables as complex variables until we finish the proof of \rt{1}. 
Let us introduce
$$
\norm{f}{r}=\max\{|f(x,y)| ; \xy \in\Delta_r\},\quad \Delta_r=\{\xy \in{\Bbb C}^2; |x|\leq r, |y|\leq r\}.
$$
Assume that $f$ and $g$ are \hol\ on $\Delta_r$. Denote
$$
B_0=\max\{\norm{f_1}{r},\norm{g}{r}\}.
$$
Assume also that
\be{3.6}
A_0=\norm{\hat f_0}{r}\leq \frac{r}{4},\qquad 0<r<1,\quad 0<\theta<1/4.  
\end{equation}

The Cauchy inequalities give
$$
|f_{\alpha,\beta}|\leq \frac{\norm{f_1}{r}}{r^{\alpha+\beta}},\qquad
|g_{\alpha,\beta}|\leq \frac{\norm{g_1}{r}}{r^{\alpha+\beta}},\quad \beta\geq 1.
$$
Hence, we get
\be{3.7}
\norm{b}{{(1-\theta)r}}\leq \sum_{k=2}^\infty\frac{B_0}{r^k}((1-\theta)r)^k\leq
\frac{B_0}{\theta}.
\end{equation}
We now have
\begin{equation*}
 \begin{split}
\norm{\hat u}{{(1-2\theta)r}}
&\leq\sum_{k\geq 2} (k+1)\mu_k((1-2\theta)r)^k\leq 2r\sum_{k\geq 2} k\mu_k((1-2\theta)r)^{k-1}\\
&\leq 2r\norm{\mu^\prime}{{(1-2\theta)r}}\leq\frac{2\norm{\mu}{{(1-\theta)r}}}{\theta},
\end{split}
\end{equation*}
in which the last inequality comes from the Cauchy formula. Now
\re{3.5}--\re{3.7} yield
\be{3.8}
\norm{\hat u}{{(1-2\theta)r}}\leq \frac{c_1B_0}{\theta^2},
\end{equation}
where, and also in the rest of discussion, $c_j>1$ stands for a 
constant. 
 With a similar computation, one can also obtain  the following  estimate 
\be{3.9}
\norm{\hat v}{{(1-2\theta)r}}\leq \frac{c_2B_0}{\theta^2}.
\end{equation}

We are ready to prove the following.
\bl{3.1}
 Let $\varphi$ be as in \rl{4.1}$.$
 Suppose that$,$ for $c_3=\max\{c_1,c_2\},$
\be{3.10}
A_0\leq \frac{r}{4},\qquad
B_0\leq \frac{\theta^3r}{4c_3}.
\end{equation}
Then we have
\be{3.11}
\varphi\colon \Delta_{(1-2\theta)r}\to\Delta_{(1-\theta)r},\qquad
\varphi^{-1}\colon \Delta_{(1-4\theta)r}\to \Delta_{(1-3\theta)r}.
\end{equation}
\el
\begin{pf} From \re{3.8}--\re{3.10}, it easy to see that $\varphi\colon
\Delta_{(1-2\theta)r}\to\Delta_{(1-\theta)r}$. 
To show the existence of the inverse mapping, we fix $(x^\prime,y^\prime)\in\Delta_{(1-4\theta)r}$ and consider the mapping
$$
T(x,y)=\left(x^\prime-u\xy ,y^\prime-v\xy \right).
$$
It is clear that $T$ maps $\nd{3}$ into itself. From \re{3.8}\ and \re{3.10}, we get
$$
\norm{u_x}{{(1-3\theta)r}} \leq \frac{\norm{u}{{(1-2\theta)r}}}{\theta r}<\frac{1}{4}.
$$
Similarly, we can verify that $\norm{u_y}{{\nnd{3}}}, \norm{v_x}{{\nnd{3}}}$ and $\norm{v_y}{{\nnd{3}}}$ are less than $1/4$. This implies that with the norm $$
\|\xy \|=\max\{|x|,|y|\}, 
$$
$T$ is a contraction mapping. By the fixed-point theorem, $T$ has a unique
fixed point $\xy $ in $\nd{3}$, which is clearly $\varphi^{-1}(x^\prime,y^\prime)$.
% The proof of \rl{3.1}\ is complete.
\end{pf}

Let us keep the notations and assumptions given in \rl{3.1}. 
Fix $(x^\prime,y^\prime)\in\nd{4}$. Then $\xy=\varphi^{-1}\xyp\in\nd{3}$.
From (\ref{eq:3.12}),
we have
\ben\begin{split}
|f(x^\prime,y^\prime)-f\xy |
&=|\int_0^1\frac{d}{dt}f(x^\prime-tu\xy ,y^\prime-tv\xy )\, dt|\\
&\leq \norm{f_x}{{(1-3\theta)r}}\norm{u}{{(1-3\theta)r}}+
\norm{f_y}{{(1-3\theta)r}}\norm{v}{{(1-3\theta)r}}.
\end{split}
\end{equation*}
From \re{3.10}, it follows that $\norm{f}{r}<r$. Now the Cauchy formula
gives
$$
\norm{f_x}{{(1-3\theta)r}}\leq\frac{1}{{3\theta }}.
$$
A similar \es\ also holds for $f_y$.  From \re{3.8}\ and \re{3.9}, we
now get
\be{ad1}%$$
|f\xyp -f\xy |\leq \frac{(c_1+c_2)B_0}{\theta^3}.
\end{equation}%$$
From \re{3.10}, we have
$$
|y+f\xy|\leq (1-3\theta)r+A_0+B_0<3r.
$$
Using \re{3.8}\ and the Cauchy formula, we get
$$
\|u_x\|_{(1-3\theta)r}\leq\frac{c_1B_0}{\theta^3r},\qquad  
\|u_y\|_{(1-3\theta)r} \leq\frac{c_1B_0}{\theta^3 r}.
$$
Hence
\be{ad2}%$$
|(y+f\xy)u_x\xy |\leq\frac{3c_1B_0}{\theta^3}.
\end{equation}%$$
 We also have
\be{ad3}%$$
 |g\xy  u_y\xy |\leq B_0\frac{c_1B_0}{\theta^3r}\leq B_0,
\end{equation}%$$
in which the last inequality is obtained by using \re{3.10}  to get rid of
one of two $B_0$'s. 
 
Substituting \re{ad1}, \re{ad2}\ and \re{ad3}\ into
\re{3.12},  we get
$$
\norm{p-f}{{(1-4\theta)r}}\leq\frac{c_4B_0}{\theta^3}.
$$
In particular, we have the estimates
\be{3.14}
\norm{p_0-f_0}{{(1-4\theta)r}}\leq\frac{c_4B_0}{\theta^3},\qquad
\norm{p_1-f_1}{{(1-4\theta)r}}\leq \frac{2c_4B_0}{\theta^3}.
\end{equation}
From \re{3.13}, one  also gets
\be{3.15}
\norm{q}{{(1-4\theta)r}}\leq\frac{c_5B_0}{\theta^3}.
\end{equation}

We are ready to prove the following.
\bp{3.2}
Let $f,g,p,q$ be as in \rl{4.1}.
Then there exist two positive constants $c_0$ and $\epsilon_0$ satisfying the following property$:$ If 
\be{3.16}
A_0=\norm{\hat f_0}{r}\leq r/4,\qquad B_0=\max\{\norm{f_1}{r},\norm{g}{r}\}\leq
\epsilon_0\theta^3r,
\end{equation}
then $\varphi$, defined by $\re{2.1},$ \tr s $\re{1}$ into $\re{n1}$ such that
\begin{gather}
A_1=\norm{\hat p_0}{{(1-5\theta)r}}\leq A_0+\frac{c_0}{\theta^4}B_0, \label{eq:3.17}\\
B_1=\max\{\norm{p_1}{{(1-5\theta)r}},\norm{q}{{(1-5\theta)r}}\}
\leq\frac{c_0}{\theta^3}B_0(1-\theta)^{d_1},\label{eq:3.18} 
\end{gather}
in which $d_1=\max\{ \ord  p_1,\ord  q\}$.
\ep
\begin{pf} We choose $\epsilon_0=1/(4c_3)$. Then \re{3.16}\ implies that   \re{3.10}\ holds. From \re{3.14}\ and Cauchy inequalities, we get
$$
|p_{k,0}-f_{k,0}|\leq\frac{c_4B_0}{\theta^3}\cdot\frac{1}
{\left((1-4\theta)r\right)^k}.
$$
Hence
$$
\norm{\widehat{p_0-f_0}}{{(1-5\theta)}}\leq\frac{c_4B_0}{\theta^3}
\sum_{k=2}^\infty\left(\frac{1-5\theta}{1-4\theta}\right)^k
\leq \frac{c_4B_0}{\theta^4 }.
$$
This gives us \re{3.17}, if we choose 
$c_0\geq c_4$. 

From \re{3.14}, we have
$$
\norm{p_1}{{(1-4\theta)r}}\leq \norm{f_1}{{(1-4\theta)r}}+\norm{p_1-f_1}{{(1-4\theta)r}}
\leq\frac{3c_4B_0}{\theta^3 }.
$$
Now, the Schwarz lemma yields
$$
\norm{p_1}{{(1-5\theta)r}}\leq \norm{p_1}{{(1-4\theta)r}}\left(\frac{1-5\theta}{1-4\theta}\right)^{d_1}.
$$
Notice that $(1-5\theta)/(1-4\theta)<1-\theta$. Therefore, we get
$$
\norm{p_1}{{(1-5\theta)r}}\leq \frac{3c_4B_0}{\theta^3}(1-\theta)^{d_1}.
$$
From \re{3.15}\ and the Schwarz inequality, we also have
$$
\norm{q}{{(1-5\theta)r}}\leq \frac{c_5B_0}{\theta^3 }(1-\theta)^{d_1}.
$$
Let us put $c_0=\max\{3c_4, c_5\}$. Then the last two inequalities yield \re{3.18}.
%The proof  of \rp{3.2}\ is complete.
%Before we end this section,  we should point out w
\end{pf}

\section{A KAM argument}\setcounter{equation}{0}
\label{section:4}
  
In this section, we shall first construct a \seq\ of \fo\ \trn s $\Phi_n$
\st\ \re{1}\ is \tr ed into \re{2}\ under the limit \trn\ of $\{\Phi_n\}$.
We shall use a KAM argument to show the convergence of the
sequence $\Phi_n$.

  We shall construct a sequence of systems
 $$\frac{dx}{dt}=y+p_n\xy ,\qquad \frac{dy}{dt}=q_n\xy ,$$
where $p_0\xy =f\xy $ and $q_0\xy =g\xy $ give the initial system \re{1}. Recursively, $p_n,q_n$ are obtained through the  \trn\ $\varphi_n$ constructed  through the approximate equations  in section 3.  

Let us decompose
$$p_n\xy =p_{n;0}\x+p_{n;1}\xy ,\qquad p_{n;0}\x=p_n(x,0).$$
Put $d_n=\min\{\ord  p_{n;1},\ord q_n\}$. We have $d_0\geq 2$.  From \rl{4.1}, we know that 
$$d_n\geq 2^n+1.$$

We now put
$$
r_n=\frac{1}{2}\left(1+\frac{1}{n+1}\right)r_0,\quad n=0,1,\ldots,
$$
in which $r_0<1$ will be determined late. Rewrite
 $$
r_{n+1}=(1-5\theta_n)r_n,\qquad \theta_n=\frac{1}{5(n+2)^2},\quad n=0,1,\ldots.
$$
We need  to choose $r_0$ so small that  all the following norms are well-defined:
$$
A_n=\norm{\hat p_{n;0}}{{r_n}},\qquad B_n=\max\{\norm{p_{n;1}}{{r_n}},\norm{q_n}{{r_n}}\}.
$$

Let us first prove a numerical result.
\bl{4.2}
Let $r_n,\theta_n,d_n$ be given as above, and let $\epsilon_0, c_0$  be
as in  \rl{3.1}. Then there exists $\epsilon_1<\epsilon_0$, which is independent of $r_0$, such that for two sequences of non-negative numbers $\{A_n^*\}_{n=0}^\infty$ and $\{B_n^*\}_{n=0}^\infty$, we have
\be{4.6}
A_n^*\leq r_n/4,\qquad B_n^*\leq\frac{\epsilon_0\theta_n^4r_n}{c_02^{n+2}},\quad n=1,2,\ldots,
\end{equation}
provided that for all $n$
\begin{gather} \label{eq:4.7}
A_{n+1}^*\leq A_n^*+\frac{c_0B_n^*}{\theta_n^4 },\qquad B_{n+1}^*\leq
\frac{c_0B_n^*}{\theta_n^3 }\left(1-\theta_n\right)^{d_{n+1}},\\
A_0^*\leq r_0/16,\qquad B_0^*\leq \epsilon_1\theta_0^4r_0.\label{eq:4.8}
\end{gather}
\end{lemma}
\begin{pf}
We put 
$$
\hat B_n=\frac{\epsilon_0\theta_n^4r_n}{c_02^{n+2}}.
$$
Clearly, we see that 
 $\hat B_{n+1}/\hat B_n\to 1$ as $n\to\infty$. On the other hand, 
$d_n>2^n$ implies that for large $n$, one has 
$d_{n+1}>\theta_n^{-2}.
$
Hence, we have
$$
(1-\theta_n)^{d_{n+1}}<(1-\theta_n)^{1/\theta_n^2}< (1 /2 )^
{ 1 /\theta_n^2},
$$
if $n$ is sufficiently large. Now, it is easy to see that
$$
\frac{c_0}{\theta_n^3 }\left(1-\theta_n\right)^{d_{n+1}}\to 0.
$$
Hence, it follows from \re{4.7}\ that there exists $n_0$ independent of the choice of $r_0$ such that
\be{4.9}
B_{n+1}^*\leq B_n^*\frac{\hat B_{n+1}}{\hat B_n},\qquad \mbox{for}\ n\geq n_0.
\end{equation}
Choose $\epsilon_1$ so small that if $ B_0^*$
satisfies the condition \re{4.8}, then
$$
B_n^*\leq \hat B_n,\qquad 0\leq n\leq n_0.
$$
Thus, \re{4.9} yields the estimate of $B_n^*$  in \re{4.6}.

As for the estimate of $A_n^*$, we have
$$
A_n^*\leq A_0^*+\sum_{j=0}^n \frac{c_0B_j^*}{\theta_j^4}.
$$ 
Using the estimate of $B_n^*$ obtained above, we get
$$
A_n^*\leq \frac{r_0}{16}+\sum_{j=0}^n\frac{\epsilon_0r_n}{2^{n+2}}.
$$
Notice that $r_0/2<r_n<r_0$. Hence
$$
\sum_{j=0}^n\frac{r_n}{2^n}<r_0\sum_{j=0}^\infty\frac{1}{2^j}=2r_0<4r_n.
$$
Thus, we obtain
$$
A_n^*\leq r_n/8+\epsilon_0r_n.
$$
We may assume that $\epsilon_0$, chosen in the proof of \rl{3.1}, is less than $1/8$. Therefore, we obtain the desired estimate of $A_n^*$.
%The proof of \rl{4.2}\ is complete.
\end{pf}
\begin{pf*}{Proof of \rt{1}}  In order to apply \rl{4.2} to $\{A_n\}_{n=0}^\infty$  and
$\{B_n\}_{n=0}^\infty$, we need to choose $r_0$. Since $\ord  f_0
\geq 2$, we may choose $r_0$ so small that  $\hat f_0$ converges for $|x|\leq r_0$,
and
\be{4.10}
\norm{\hat f_0}{{r_0}}\leq r_0/16.
\end{equation}
By choosing a smaller $r_0$ if it is
necessary, we may also assume that
$f_1,g$ are \hol\ functions on $\Delta_{r_0}$, and 
\be{4.11}
\norm{f_1}{{r_0}}\leq \epsilon_1\theta_0^4r_0,\qquad\norm{g}{{r_0}}\leq \epsilon_1\theta_0^4r_0.
\end{equation}

From \re{4.10}\ and \re{4.11}, it follows that \re{3.16}\ is satisfied. Hence, 
\rp{3.2}\ says that $\{A_0,A_1\}$ and $\{B_0,B_1\}$ satisfy \re{4.7}. One also sees 
that two initial conditions in \re{4.8} follow from \re{4.10}\ and \re{4.11}. Now, 
\rl{4.2} implies that 
$$
A_1\leq r_1/4,\qquad B_1\leq
\frac{\epsilon_0\theta_1^4r_1^2}{c_02^3}.
$$
In particular, this gives us two initial conditions  in \re{3.10}\ for the new system defined by $p_1$ and $q_1$. Hence, we may
 apply \rp{3.2}\ again. By repeating  this process, we can prove that 
\be{4.12}
A_n\leq r_n/4,\qquad B_n\leq
\frac{\epsilon_0\theta_n^4r_n}{c_02^{n+2}},\quad n=0,1,\ldots.
\end{equation}
Therefore, \rl{3.1}\ gives us
$$
\varphi_n^{-1}\colon \Delta_{r_{n+1}}\to\Delta_{r_n},\qquad n=0,1,\ldots.
$$

Notice that $r_n\geq r_0/2$. We have
$$
\Phi_n^{-1}=\varphi_0^{-1}\circ\varphi_1^{-1}\circ\ldots\circ\varphi_n^{-1}\colon 
\Delta_{\frac{1}{2}r_0}\to\Delta_{r_0}.
$$
Hence
$$
\norm{\Phi_{n+1}^{-1}-\Phi_n^{-1}}{{\frac{1}{2}r_0}}\leq 2r_0,
$$
where the norm on the left side is defined to be the maximum of  norms of
two components. On the other hand, we know from \re{3.3}\ that 
each component of
$\Phi_{n+1}^{-1}-\Phi_n^{-1}$ vanishes with order at least $d_n>2^n$. From
Schwarz Lemma,  it follows that
$$
\norm{\Phi_{n+1}^{-1}-\Phi_n^{-1}}{{\frac{1}{4}r_0}}\leq
r_0\left(\frac{1}{2}\right)^{2^n}.
$$
Therefore, the \seq\ $\Phi_n^{-1}$ converges to a \trn\
$\Phi_\infty^{-1}$.

It is clear that $\Phi_\infty$ \tr s the system \re{1}\ into
a system of the form \re{2}. Since the normalized \trn s form a group, we
see that $\Phi_\infty$ is still a normalized \trn. In section 2, we have
seen that $\Phi$ is the unique formal \trn\ which \tr s \re{1}\ into
\re{2}. Therefore, we obtain that $\Phi=\Phi_\infty$, so $\Phi$ is a convergent
\trn.
% The proof of \rt{1}\ is complete.
\end{pf*}

\section{Embeddability of parabolic mappings}\setcounter{equation}{0}
\label{section:5}

In this section, we shall investigate the relation between the embeddability
of parabolic mappings as the time-1 mappings and the convergence of
normalization for parabolic \trn s. We shall first show that a parabolic
mapping is
formally \li\ if and only if it is embeddable as a time-1 mapping of  a 
formally \li\ \vf. 

Let us put
$$
x=\left(\begin{array}{c} x_1\\ x_2\end{array}\right) , \qquad 
A=\jz{
{cc} 0 & 1\\ 0&0
}.
$$
Rewrite \re{1}\  as
\be{5.2}
\frac{dx}{dt}=Ax+\sum_{|I|\geq 2}F_Ix^I\equiv F(x),
\end{equation}
where each $F_I$ is a constant matrix of $2$ by $1$, and  $x^I=x_1^\alpha
x_2^\beta$ for $I=(\alpha,\beta)$.

We first assume that
$F\x$ is only given by formal \pows. Let $\varphi_t$ be a family of \fo\ \trn s
generated by the formal \vf\  \re{5.2}, i.~e.
\be{flow}
\frac{d}{d t}\varphi_t=F\circ\varphi_t,\qquad \varphi_0=\mbox{Id}.
\end{equation}
Since the linear part of $\varphi_t$ with respect to $x_1$ and $x_2$ is determined by the matrix $A$, then we have the expansion
$$
\varphi_t(x)=e^{At}x+\sum_{|I|\geq 2}B_I(t)x^I,\qquad B_I(0)=0,
$$
where $B_I(t)$ is a matrix of $2$ by $1$ given by formal \pows, and
$$
e^{At}=\sum_{n=0}^\infty\frac{A^n}{n!}t^n=\jz{{cc} 1&t\\ 0&1}.
$$
Now, \re{flow}\ takes the form
\be{5.3}
B_I^{\prime}(t)=AB_I(t)+b_I(t),\qquad B_I(0)=0,
\end{equation}
where $b_I(t)$ depends only on $A$ and $B_J(t)$ with $|J|<|I|$. More
precisely, we have
\be{5.4}
b_I(t)=\left\{\sum_{|J|=|I|}F_J(e^{At}x)\right\}_I+\tilde{b}_I(t),
\end{equation}
where  $\{\cdot\}_I$ denotes the matrix of  coefficients for $x^I$, and 
$$
\tilde{b}_I(t)\equiv 0,\qquad \mbox{if}\ F_J=0,\ \mbox{for all }\ |J|<|I|.
$$
The solution $B_I$ to \re{5.3}\ is given by
\be{5.5}
B_I(t)=e^{At}\int_0^te^{-At}b_I(t)\ dt.
\end{equation}
Hence, $B_I(t)$ are real \an\ functions defined on the whole real line. In particular, we see that $\varphi_t$ is a family of formal \trn s defined for $-\infty<t<\infty$.
One also notices that if $\tilde{v}=\Psi_*v$ for a formal \trn\ $\Psi$, then the 1-parameter family  of formal \trn s generated by $\tilde{v}$ are given by $\Psi\circ\varphi_t\circ\Psi^{-1}$.

We need the following lemma.
\bl{5.2}
Let $v$ be a formal \vf\ defined by $\re{5.2},$ and $\varphi_t$ the
1-parameter family of formal \trn s generated by $v$. 
Assume that  $\re{5.2}$ is not a linear system. Then $\varphi_t$ is not a linear \trn\ for all $t\neq 0$.
\el
\begin{pf}
We assume that there is $I_0=(\alpha_0,\beta_0)$ such that
$$
F_{I_0}\neq 0,\qquad F_J= 0,
$$
for $|J|<|I_0|$, or for $J=(\alpha,|I_0|-\alpha)$ with $\alpha>\alpha_0$. Then from \re{5.4}, it follows that
$$
b_{I_0}(t)=F_{I_0}\neq 0.
$$
Next, we use the formula \re{5.5} and get
$$
B_{I_0}(t)=\jz{{cc}t&t^2/2\vspace{.5ex}\\ 0&t}b_{I_0}\neq 0,\qquad \mbox{for}\ t\neq 0,
$$
which implies that $\varphi_t$ is not a linear \trn\ for each
$t\neq 0$.
%Therefore, the proof of \rl{5.2}\ is complete.
\end{pf}

We now consider a parabolic \trn\  
\be{5.1}
 \varphi\xy =T(x,y)+O(2),\qquad T\xy=(x+y,y).
\end{equation}
In~\cite{gong1}, it was proved that there exist real \an\ parabolic
mappings which are not \li\ by any convergent transformation. 
In fact, the parabolic mappings are constructed through 
 a pair of real analytic glancing hypersurfaces. On the other hand,
Melrose~\cite{melrose1} showed that a pair of smooth glancing hypersurfaces
can always be put into a certain normal form by smooth transformations; and consequently, the parabolic
mappings coming from a pair of smooth glancing hypersurfaces are always \li\ by
smooth \trn s.
 Therefore, we can state the following.
\bt{5.1}
There exists a smoothly linearizable  real \an\ \trn\ $\varphi$ of the form $\re{5.1}$, which cannot be \tr ed into $T$ by any convergent \trn.
\et
Now, we see that \nrc{3}\ follows from the following.
\bp{5.3}
Let $\varphi$ be a real \an\ \trn\ of the form $\re{5.1}.$ Assume that
$\varphi$ is formally equivalent to $T$. Then $\varphi$ is a time-1
mapping of a real \an\ \vf\ of the form $\re{5.2}$, if and only if $\varphi$ can be \tr ed into $T$ through a convergent \trn.
\ep
\begin{pf}
Obviously, $\varphi$ is   embeddable  if it is linearizable through convergent \trn s. We now assume that
$\varphi=\varphi_1$ for a 1-parameter family of  \trn s $\varphi_t$
generated by a real \an\ $v$ of the form \re{5.2}. Let $\Phi$ be a formal
\trn\ which linearizes $\varphi$. This implies that  the time-1 mapping
of the formal \vf\ $\Phi_*v$ is a linear \trn. From \rl{5.2}, it follows 
 that  $\Phi_*v$ is a linear \vf.
Now, \rt{1}\ implies that $v$ is \li\ by a convergent
\trn\ $\phi$. Using \rl{5.2} again, we know that $\varphi$ is also \li\ by
the same \trn\ $\phi$. 
%Therefore, the proof of \rp{5.3} is complete.
\end{pf}
\bibliographystyle{plain}

\end{document}